\documentclass[a4paper,10pt]{amsproc}

\usepackage{amsfonts, amsmath, amssymb, graphicx, mathrsfs}
\usepackage[all]{xy}
\usepackage{enumerate}
\newdir{ >}{!/6pt/@{ }*:(1,0.2)@_{>}*:(1,-0.2)@^{>}}

\theoremstyle{plain}

\newtheorem{thm}{Theorem}[section]
\newtheorem{lem}[thm]{Lemma}
\newtheorem{proposition}[thm]{Proposition}
\theoremstyle{definition}
\newtheorem{defn}[thm]{Definition}

\newtheorem{cor}[thm]{Corollary}

\newtheorem{counter example}[thm]{Counter Example}

\newtheorem{exam}[thm]{Example}

\numberwithin{equation}{section}

\author[S. Mandal]{Samir Ch Mandal}
\address{Department of Pure Mathematics, University of Calcutta, 35, Ballygunge Circular Road, Kolkata 700019, West Bengal, India}
\email{samirchmandal@gmail.com}

\author[S. Bag]{Sagarmoy Bag}
\address{Department of Mathematics, Bangabasi Evening College, 19 Rajkumar Chakraborty Sarani, Kolkata 700009, West Bengal, India}
\email{sagarmoy.bag01@gmail.com}

\author[D. Mandal]{Dhananjoy Mandal}
\address{Department of Pure Mathematics, University of Calcutta, 35, Ballygunge Circular Road, Kolkata 700019, West Bengal, India}
\email{dmandal.cu@gmail.com}

\thanks{ The first author thanks to UGC, India, for financial support.\\
  Corresponding author.}


\begin{document}

\title[ A note on the ring $C(X)_F$]{A note on the rings of functions which are discontinuous on some finite sets }

%    General info
%%%%%%%%%%%%%%%%%%%%%%%%%%%%%%%%%%%%%%%%%%%%%%%%%%%

%                                                                                                                           %
%         Please use the current 2010 Mathematics Subject Classification:             %
%         http://www.ams.org/mathscinet/msc/                                                        %
%         http://www.zentralblatt-math.org/msc/en/                                                 %
%%%%%%%%%%%%%%%%%%%%%%%%%%%%%%%%%%%%%%%%%%%%%%%%%%%

\thanks {}
\maketitle	

\begin{abstract}
In this paper, we study some properties of the ring $C(X)_F$ of all real valued functions which are continuous except on some finite subsets of $X$. We show that $C(X)_F$ is closed under uniform limit if and only if the set of all non-isolated points of $X$ is finite. We also initiate and investigate the zero divisor graph of the ring $C(X)_F$.
\end{abstract}	
\section{Introduction}
All our spaces are $T_1$ unless otherwise  stated. Let $C(X)$ be the ring of all real valued continuous functions on a topological space $X$ and it has been investigated extensively [see \cite{gill}]. The ring $C(X)_F$ of all real valued functions on $X$ which are discontinuous on some finite subset of $X$ was first introduced  and studied in \cite{main} and very recently some interesting properties of $C(X)_F$ have been investigated also by  M. R. Ahmadi Zand and Z. Khosravi \cite{appl}. The ring $B_1(X)$ of all real valued Baire class one functions on $X$ was introduced and investigated by A. Deb Ray and A. Mondal \cite{bare1,bare2}. Our intention in the present article is to pursue research on the ring $C(X)_F$.  
\par In  Section 2, we study some properties of $C(X)_F$ and examine the relation between the rings $C(X)_F$ and $B_1(X)$. In this section, we also establish that $C(X)_F$ is closed under uniform limit if and only if the set of all non-isolated points of $X$ is finite. 
\par Another interesting ring $T'(X)$ of all real valued functions which are continuous on some open dense subset of $X$ was extensively studied in \cite{blum}. We have investigated as to when $T'(X)$ is closed under uniform limit in Section 2 of this article. 
\par For a commutative ring $R$ with unity, let $\Gamma(R)$ be the set of all nonzero  zero divisors of $R$. Consider an edge between two vertices $a$, $b$ in $\Gamma(R)$ $($ i.e., $a$ and $b$ are adjacent$)$ if $ab=0$. $\Gamma(R)$ is called zero divisor graph of the ring $R$ \cite{fdapsl}. For $a,b\in \Gamma(R), \:d(a,b)$, is defined to be the length of the smallest path between $a$ and $b$. Diameter of $\Gamma(R)$, denoted by $diam(\Gamma(R))$ is defined by the maximum of all possible $d(a,b),$ for $a,b\in\Gamma(R)$. A graph $\Gamma(R)$ of a commutative ring $R$ with unity is called triangulated if each vertex of $\Gamma(R)$ belongs to a cycle of length 3. The zero divisor graph $\Gamma (C(X))$ of $C(X)$ was introduced by F. Azarpanah and M. Motemedi \cite{zero}. They have observed the inter relations among the ring properties of $C(X)$, the graph properties of $\Gamma (C(X))$ and the topological properties of $X$.
\par In the third section, we study some interesting results on the zero divisor graph of the ring $C(X)_F$. We have shown that $diam(\Gamma(C(X)_F))$ is 3 if $X$ contains at least 3 points [See Corollary \ref{diam}].  We have established that $\Gamma(C(X)_F)$ is never triangulated [Theorem \ref{try}].

 We conclude our paper in Section 4 by pointing out the similarities and dissimilarities of the rings $C(X)$ and $C(X)_F$, based on the results obtained in this paper.

 Throughout this article, the set of real numbers, rational numbers and natural numbers are denoted by $\mathbb{R}$, $\mathbb{Q}$ and $\mathbb{N}$ respectively.
 
\section{Some properties of  $C(X)_F$}
In this section, we study some properties of $C(X)_F$ and investigate relation between the two rings $C(X)_F$ and $B_1(X)$.  \\
 
Let us give some preliminary definitions and results which will be used subsequently in this section. 
\begin{defn}%[Perfectly Normal Space]
	A topological space $X$ is called perfectly normal if it is normal and every closed subset of $X$ is a $G_\delta$-set. 
\end{defn}
\begin{defn}%[$P$-space]
	A topological space $X$ is called a $P$-space if every $G_\delta$-set  is open in $X$.
\end{defn}
\begin{defn}%[Baire Space]
	A topological space $X$ is said to be Baire space if any countable intersection of open dense subsets in $X$ is dense in $X$.
\end{defn}
The following theorem \cite{vesel} gives a characterization of Baire one functions. 
\begin{thm}[\cite{vesel}]\label{bare-condition}
	For a normal topological space $X, B_1(X) = F_\sigma(X)$, where $F_\sigma(X)=\{f\in \mathbb{R}^X: f^{-1}(U) \text{ is an }F_\sigma\text{-set for any open set } U\subseteq \mathbb{R}\}$. 
\end{thm} 
\begin{thm}[Lemma 5.3, \cite{bare1}] \label{bare-uniform}
	Suppose $\{f_n\}$ is a sequence of functions in $B_1(X)$ which converges uniformly to $f:X\rightarrow \mathbb{R}$. Then $f\in B_1(X)$. In  other words, $B_1(X)$ is closed under uniform limit.
\end{thm}
Now we are in a position to describe our result. 
\begin{thm}\label{perfect}
	If $X$ is a perfectly normal space, then $C(X)_F\subseteq B_1(X)$.
\end{thm}
\begin{proof}
	Let $f\in C(X)_F.$ Then there exists a finite subset $F$ of $X$ such that $f|_{X\setminus F}\in C(X\setminus F)$. Now take any open set $U\subseteq \mathbb{R}$. Then $f^{-1}(U)=(f|_{X\setminus F})^{-1}(U)\cup K, $ for some finite set $K\subseteq F.$  Since $(f|_{X\setminus F})^{-1}(U)$ is open in $X\setminus F$ and $X\setminus F$ is open in $X$, it follows that $(f|_{X\setminus F})^{-1}(U)$ is open in $X$. Now as $X$ is perfectly normal and $(f|_{X\setminus F})^{-1}(U)$ is open in $X$,  it is an $F_\sigma$-set. On the other hand, $K$ being a finite set, is an $F_\sigma$-set. Since the union of two $F_\sigma$-sets is an $F_\sigma$-set,  $f^{-1}(U)$ is an $F_\sigma$-set. As $X$ is perfectly normal, it is normal too. Thus we can apply Theorem \ref{bare-condition} to see that $f\in B_1(X)$. Hence $C(X)_F\subseteq B_1(X)$.
\end{proof}
For any normal space $X$, $C(X)_F$ may  not be contained in $B_1(X)$, as is seen in the following example:

\begin{exam} Consider the set $\mathbb{R}$ with discrete topology. Let $\mathbb{R}^*=\mathbb{R}\cup \{\infty\}$ be the one-point compactification of $\mathbb{R}$, where $\infty\notin \mathbb{R}$. Then $\mathbb{R}^*$ is normal but not perfectly normal. Indeed the set $\{\infty\}$ is not a $G_\delta$-set. This implies that $\mathbb{R}$ is not an $F_\sigma$-set.\\ Define $f\in C(\mathbb{R}^*)_F$ as $f=\chi_\mathbb{R}$, where $\chi_\mathbb{R}$ is the characteristic function on $\mathbb{R}$. Then $f^{-1}(0.5,1.5)=\mathbb{R}$ is not an $F_\sigma$-set. Thus appealing to Theorem \ref{bare-condition}, the function $f \notin B_1(\mathbb{R}^*)$.
\end{exam}
In the following, we produce an example that reflects that $B_1(X)$ strictly contains $C(X)_F$.
\begin{exam}
Let  $X=[0,1]$ be the subspace of the real line $(\mathbb{R}, \tau_u)$, where $\tau_u$ is the Euclidean topology on $\mathbb{R}$. Consider $[0,1]\cap \mathbb{Q}=\{r_1,r_2,\cdots \}$. For $i\in \mathbb{N},$ write $r_i=\frac{p_i}{q_i}$, where $p_i,q_i\in \mathbb{Z}$ and $gcd(p_i,q_i)=1$.\\
  Now we define a sequence of functions $\{f_n\}$  on $X$ as follows:
	\begin{align*}
		f_n(x)&=\left\{\begin{array}{cl}
			\frac{1}{q_i},& \text{if }x=r_i=\frac{p_i}{q_i},\: i=1,2,\cdots,n\\
			0, & \text{elsewhere.}
		\end{array}\right.
	\end{align*}
	 Then each $f_n\in C(X)_F$. We claim that the sequence $\{f_n\}$ converges uniformly to $f$ on  $X$, where $f$ is as follows:
	 	\begin{align*}	 
	 	f(x)&=\left\{\begin{array}{cl}
	 		\frac{1}{q},& \text{if }x=\frac{p}{q},\: p,q\in \mathbb{Z} \text{ and } gcd(p,q)=1\\
	 		0, & \text{elsewhere.}
	 	\end{array}\right.
	 \end{align*} 
In fact, let $\epsilon>0.$ We need to find some $k\in \mathbb{N}$ such that 
		\[|f_n(x)-f(x)|<\epsilon \text{ for all } x\in [0,1] \text{ and for all }n\geq k.\]
		Now there are finitely many points of the form $\frac{p}{q}\in [0,1]$ such that $\frac{1}{q}\geq \epsilon.$ Let these points lie in $\{r_1,r_2,\cdots,r_{k-1}\}$. Again for all $n\geq k, $ if $ r_n=\frac{p}{q}$, then $\frac{1}{q}<\epsilon$. So for all $n\geq k,$ we have 
		\[|f_n(x)-f(x)|=\left\{\begin{array}{cl}
			\frac{1}{q}, & \text{if }x=\frac{p}{q}, x\in \{r_{n+1},r_{n+2},\cdots\}\\
			0, & \text{elsewhere.}
		\end{array}\right.\]
		This shows that  $|f_n(x)-f(x)|<\epsilon \text{ for all } x\in [0,1] \text{ and for all }n\geq k.$ \\
	Clearly $f\notin C(X)_F$ as $f$ is discontinuous on the whole of $\mathbb{Q}\cap[0,1]$ which is an  infinite set in $[0,1].$  Again $X$ being perfectly normal, each $f_n\in B_1(X)$ [see Theorem  \ref{perfect}] and hence by Theorem \ref{bare-uniform}, $f\in B_1(X)$. Thus $f\in B_1(X)\setminus C(X)_F$.

\end{exam}
 The above example shows that $C(X)_F$ is not always closed under uniform limit. In the following result, we have furnished a necessary and sufficient condition on a topological space $X$ so that $C(X)_F$ is closed under the uniform limit. 
 
\begin{thm}\label{uniform}
	Any uniformly convergent sequence of functions in $C(X)_F$ has the limit function in $C(X)_F$ if and only if the number of non-isolated points in $X$ is finite.
	
\end{thm}
\begin{proof}
	First we will show that if the number of non-isolated points is infinite then there exists a sequence $\{f_n\}$ of functions in $C(X)_F$ which does not converge uniformly  to a function in $C(X)_F$. For this, let $A=\{a_1, a_2, \cdots\}$ be an infinite set of non-isolated points in $X$. 
	\par Define \[f_n(x)=\left\{\begin{array}{lcl}
		\frac{1}{k},&& \text{if }x=a_k\:,k=1,2,\cdots,n-1\\
		0, && \text{elsewhere. } 
	\end{array}\right.\]
Then each $f_n \in C(X)_F$, $n\in \mathbb{N}$. We now show that $\{f_n\}$  converges uniformly to $f$, where $f$ is defined by 
\[f(x)=\left\{\begin{array}{lcl}
	\frac{1}{k},&& \text{if }x=a_k\in A\\
	0, && \text{if } x\notin A.
\end{array}\right.\] 
Now we have \[
|f_n(x)-f(x)|=\left\{\begin{array}{lcl}
	\frac{1}{k},&&\text{if }x=a_k,\: k\geq n\\
	0,&& \text{elsewhere.} 
\end{array}\right.\]
Choose any $\epsilon>0.$ Then we can find a natural number $k$ such that $\frac{1}{k}<\epsilon.$
So  $|f_n(x)-f(x)|\leq \frac{1}{n}\leq \frac{1}{k}<\epsilon$ for all $x\in X$ and for all $n\geq k.$ Hence $\{f_n\}$ converges uniformly to $f$. We note that each $f_n \in C(X)_F$ whereas $f\notin C(X)_F$ as the set of discontinuity of $f$ is at least the infinite set  $A$. Thus $\{f_n\}$ does not converge uniformly to a function in $C(X)_F$.
\par Conversely, let $K=\{b_1,b_2,\cdots,b_k\}$ be the set of all non-isolated points of $X$. Then $C(X)_F=\mathbb{R}^X$. Certainly then $C(X)_F$ is closed under uniform limit. This concludes the proof.   
\end{proof}
Combining Theorem \ref{uniform} and Theorem 3.4 \cite{main}, we have the following result.
\begin{thm}
	For a topological space $X$, the following are equivalent.\\
	$(i)$ $X$ has finite number of non-isolated points.\\
	$(ii)$ $C(X)_F$ is closed under uniform limit.\\
	$(iii)$ There is a topological space $Y$ such that $C(X)_F \cong C(Y)$.
\end{thm}

\par In this context we are interested to study as to when $T'(X)$ is closed under uniform limit. 
\begin{thm}\label{tdash1}
	Let $X$ be a Baire space, $P$-space and $\{f_n\}$ be a sequence of functions in $T'(X)$ that converges uniformly to $f$ on $X$. Then $f\in T'(X)$.
\end{thm}
\begin{proof}
	Let $D_n$ be an open dense subset of $X$ such that $f_n|_{D_n}$ is continuous. Then $D=\cap D_n$ is a dense subset of $X$, as $X$ is a Baire space. Also $X$ is a $P$-space so  $D$ is an open set. Thus $D$ becomes an open dense subset of $X$. Now from the hypothesis $f_n|_D$ converges uniformly to $f|_D$, where each $f_n|_D$ is continuous, as $D\subseteq D_n$ for all $n\in \mathbb{N}$. Therefore $f|_D$ is a continuous function. Therefore $f\in T'(X)$. This proves that $T'(X)$ is closed under uniform limit under the given condition.
\end{proof}
Below we furnish an example of a non-discrete space $X$ which is both Baire and $P$-space to show that our result can be applied to some nontrivial topological space. For this let us take $X=\mathbb{R}$ and define a topology $\tau $ on $X$ in the following way. A subset $U$ of $X$ is open in $X$ if it does not contain 0 or if it contains zero then $X\setminus U$ must be countable set. Then by Problem 4N.1 (\cite{gill}),  $(X,\tau)$ is non-discrete $P$-space.  Only open dense subsets of $X$ are $X$ and $X\setminus\{0\}$. This shows that $(X,\tau)$ is a Baire space too. 

\par  In \cite{appl}  it is established that $C(X)_F = T'(X)$  under some conditions on $X$. For quick reference it is mentioned below. 
\begin{thm}[Theorem 3.5, \cite{appl}]\label{tdash2}
	Let $X$ be an arbitrary topological space and let $I(X)$ denote the set of all isolated points of $X$. Then $T'(X)=C(X)_F$ if and only if for  every point $x$ of $X\setminus I(X),\{x\}$ is nowhere dense and for every open dense subset $U$ of $X$ we have, $X\setminus U$ is finite. 
\end{thm} 
In the set up when $X$ is a $P$-space as well as  a Baire space, we have the following  nice characterization of the space $X$ for $T'(X)=C(X)_F$. 
\begin{thm}\label{tdash3}
	Let $X$ be a Baire space, $P$-space. Then $T'(X)=C(X)_F$ if and only if the number of non-isolated points in $X$ is finite.
\end{thm}
\begin{proof}
	Let $T'(X)=C(X)_F$. Then as $T'(X)$ is closed under uniform limit, $C(X)_F$ is also so. So by Theorem \ref{uniform}, $X$ has only finitely many non-isolated points. 
	\par Conversely, let the number of non-isolated points in $X$ be finite. Then $C(X)_F=\mathbb{R}^X$. We claim that $T'(X)=\mathbb{R}^X$. For this let $f\in \mathbb{R}^X$. Then the complement of the set of non-isolated points is open dense in $X$. Clearly $f$ is continuous on this set. Hence $f\in T'(X)$. This completes the proof.  
\end{proof}

\newpage
\section{Zero Divisor Graph on $C(X)_F$}
Consider the graph $\Gamma(C(X)_F)$ of the ring $C(X)_F$, where the set of all vertices $V$ is the collection of  all nonzero zero divisors in the ring $C(X)_F$ and there is an edge between two vertices $f,g $ $($ i.e., $f$ and $g$ are adjacent) if $fg=0$. For $f\in C(X)_F$, $\mathcal{Z}(f)=\{x\in X: f(x)=0\}$ is called the zero set of $f$. For a $T_1$ space  $X$ and $x\in X$, the characteristic function $\chi_{\{x\}}$ defined by $\chi_{\{x\}}(y)=0$ if $y\neq x$ and $\chi_{\{x\}}(x)=1$ is a member of $C(X)_F$.
\par The following result gives a characterization of the existence of edge between two vertices of $\Gamma(C(X)_F)$ in terms of their zero sets.

\begin{lem}\label{deg1}
In the graph $\Gamma(C(X)_F)$, for any two vertices $f,g$ there is an edge between $f$ and $g$ if and only if $\mathcal{Z}(f)\cup\mathcal{Z}(g)=X$.
\end{lem}
\begin{proof}
	First we assume that $\mathcal{Z}(f)\cup\mathcal{Z}(g)=X$. Then clearly  $fg=0$. So there is an edge between $f$ and $g$.
	\par Conversely, let there be an edge between $f$ and $g$. Then  $fg=0$. This implies that  $\mathcal{Z}(fg)=\mathcal{Z}(f)\cup\mathcal{Z}(g)=X$.
\end{proof}
\begin{lem} \label{deg 2}
In the	graph $\Gamma(C(X)_F)$, for any two vertices $f,g$ there is another vertex $h\in \Gamma(C(X)_F)$ such that $h$ is adjacent to both $f$ and $g$ if and only if $\mathcal{Z}(f)\cap\mathcal{Z}(g)\neq \emptyset$. 
\end{lem}
\begin{proof}
First we assume that there is a vertex $h$ such that $h$ is adjacent to both $f$ and $g$. Then $hf=0$ and $hg=0$.	As $h$ is non-zero, there is a point $x_0\in X$ such that $h(x_0)\neq 0$. Then obviously, $f(x_0)=0$ and $g(x_0)=0$. So $x_0\in \mathcal{Z}(f)\cap\mathcal{Z}(g)$. Hence $\mathcal{Z}(f)\cap\mathcal{Z}(g)\neq \emptyset$. 
\par Conversely, assume that $\mathcal{Z}(f)\cap\mathcal{Z}(g)\neq \emptyset$. Now let $y\in \mathcal{Z}(f)\cap\mathcal{Z}(g)$. Consider $h=\chi_{\{y\}}$. Then clearly $h\in C(X)_F$ and $hf=0$ and $hg=0$. So $h$ is adjacent to both $f$ and $g$.
\end{proof} 
\begin{lem} \label{deg 3}
In the	graph $\Gamma(C(X)_F)$,	for any two vertices $f,g$ there are distinct vertices $h_1$ and $h_2$ in $\Gamma(C(X)_F)$ such that $f$ is adjacent to $h_1$, $h_1$ is adjacent to $h_2$ and $h_2$ is adjacent to $g$ if $\mathcal{Z}(f)\cap\mathcal{Z}(g)= \emptyset$.
\end{lem}
\begin{proof}
	As $\mathcal{Z}(f)\cap\mathcal{Z}(g)= \emptyset$, we can choose distinct points $x,y\in X$ such that $x\in \mathcal{Z}(f)$ and $y\in \mathcal{Z}(g)$. Consider two functions $h_1=\chi_{\{x\}}$ and $h_2=\chi_{\{y\}}$. Then we have $\mathcal{Z}(h_1)=X\setminus \{x\}$ and $\mathcal{Z}(h_2)=X\setminus \{y\}.$ So we get, $\mathcal{Z}(h_1)\cup \mathcal{Z}(f)=X$, $\mathcal{Z}(h_2)\cup \mathcal{Z}(g)=X$ and $\mathcal{Z}(h_1)\cup \mathcal{Z}(h_2)=X$. Therefore by  Lemma \ref{deg1}, we can say that $f$ is adjacent to $h_1$, $h_1$ is adjacent to $h_2$ and $h_2$ is adjacent to $g$.
\end{proof}
\begin{defn}
	For two vertices $f,g$ in $\Gamma(C(X)_F)$ , $d(f,g)$ is defined as the length of the smallest path between $f$ and $g$.
\end{defn}
\begin{thm}\label{deg-thm}
	For any two vertices $f,g$ in the	graph $\Gamma(C(X)_F)$, we have the following results: 
	\begin{enumerate}[$(i)$]
		\item $d(f,g)=1$ if and only if $\mathcal{Z}(f)\cup\mathcal{Z}(g)=X$.
		\item $d(f,g)=2$ if and only if $\mathcal{Z}(f)\cup\mathcal{Z}(g)\neq X$ and $\mathcal{Z}(f)\cap\mathcal{Z}(g)\neq \emptyset$.
		\item $d(f,g)=3$ if and only if $\mathcal{Z}(f)\cup\mathcal{Z}(g)\neq X$ and $\mathcal{Z}(f)\cap\mathcal{Z}(g) = \emptyset$.
	\end{enumerate}
\end{thm}
\begin{proof} (i) It follows from Lemma \ref{deg1}.
	\par  (ii) Let $d(f,g)=2.$ Then $f$ and $g$ are not adjacent to each other. So by Lemma \ref{deg1}, we have $\mathcal{Z}(f)\cup\mathcal{Z}(g)\neq X$. On the other hand, there is a vertex $h\in \Gamma(C(X)_F)$ such that $h$ is adjacent to both $f$ and $g$. So by Lemma \ref{deg 2}, we have $\mathcal{Z}(f)\cap\mathcal{Z}(g)\neq \emptyset$.
	\par Conversely, let $\mathcal{Z}(f)\cup\mathcal{Z}(g)\neq X$ and $\mathcal{Z}(f)\cap\mathcal{Z}(g)\neq \emptyset$. Then by Lemma \ref{deg1}, there is no edge between $f$ and $g$. Also by Lemma \ref{deg 2}, there is a vertex $h$ which is adjacent to both $f$ and $g$. So $d(f,g)=2.$
	\par  (iii) Let $d(f,g)=3$. Then by Lemma \ref{deg1} and Lemma \ref{deg 2}, we get  $\mathcal{Z}(f)\cup\mathcal{Z}(g)\neq X$ and $\mathcal{Z}(f)\cap\mathcal{Z}(g) = \emptyset$.
	\par Conversely, let $\mathcal{Z}(f)\cup\mathcal{Z}(g)\neq X$ and $\mathcal{Z}(f)\cap\mathcal{Z}(g) = \emptyset$. Then by Lemma \ref{deg1} and Lemma \ref{deg 2}, we know that $f$ and $g$ are not adjacent to each other and there is no common vertex $h$ which is adjacent to both $f$ and $g$. So $d(f,g)\geq 3.$ As $\mathcal{Z}(f)\cap\mathcal{Z}(g) = \emptyset$, applying Lemma \ref{deg 3}, there are two distinct vertices $h_1$ and $h_2$ such that $f$ is adjacent to $h_1$, $h_1$ is adjacent to $h_2$ and $h_2$ is adjacent to $g$. Hence we conclude that $d(f,g)=3.$
\end{proof}
\begin{defn}
	The diameter of the graph $\Gamma(C(X)_F)$ is denoted by $diam(\Gamma(C(X)_F))$, defined by the maximum of all possible $d(f,g)$. Here the girth of the graph $\Gamma(C(X)_F)$ is denoted by $gr(\Gamma(C(X)_F))$, defined as the length of the smallest cycle in $\Gamma(C(X)_F)$. If there does not exist any circle in the graph $\Gamma(C(X)_F)$, then we declare $gr (\Gamma(C(X)_F))=\infty.$ 
\end{defn}

\begin{thm}\label{diam}
	If $X$ contains at least three elements, then $diam(\Gamma(C(X)_F))=3$ and $gr(\Gamma(C(X)_F))=3.$
\end{thm}
\begin{proof}
	Let us take three distinct points $x,y,z$ from $X$. Consider the functions  $f=1-\chi_{\{x\}}$ and $g=1-\chi_{\{y\}}$. Then $\mathcal{Z}(f)=\{x\}$ and also $\mathcal{Z}(g)=\{y\}$. So we have, $\mathcal{Z}(f)\cup \mathcal{Z}(g)\neq X$ as $z$ does not lie in the union and $\mathcal{Z}(f)\cap\mathcal{Z}(g)=\emptyset$. So by Theorem \ref{deg-thm}(iii), $d(f,g)=3.$ Again we know that $d(f,g)\leq 3$ for all $f,g\in \Gamma(C(X)_F)$. Hence we have $diam(\Gamma(C(X)_F))=3.$
	\par For the second part, take $h_1=\chi_{\{x\}}$, $h_2=\chi_{\{y\}}$ and $h_3=\chi_
	{\{z\}}$. Then the union of any two zero sets among $\mathcal{Z}(h_1)$, $\mathcal{Z}(h_2)$ and $\mathcal{Z}(h_3)$ is $X$. Then we get a cycle of length three containing the vertices $h_1$, $h_2$ and $h_3$. As, there is no loop in the graph $\Gamma(C(X)_F)$, we have proved that the girth $gr(\Gamma(C(X)_F))=3$.
\end{proof}
\begin{thm}
	In the graph $\Gamma(C(X)_F)$, $|X|=2$ if and only if $diam(\Gamma(C(X)_F))=2$ if and only if $gr(\Gamma(C(X)_F))=\infty$. 
\end{thm}
\begin{proof}
	Let $X$ contains only two points, say $X=\{x,y\}$. Then for any vertex $f$ of $\Gamma(C(X)_F)$, we must have $\mathcal{Z}(f)$ is singleton. Then choose $f=\chi_{\{x\}}$ and $g=\chi_{\{y\}}$. Then $f$ and $2f$ are not adjacent to each other whereas $g$ is adjacent to both $f$ and $2f$. This shows that if for two functions their zero sets are same then they have distance two and if their zero sets are not same then they are adjacent to each other. Hence we can conclude that $diam(\Gamma(C(X)_F))=2.$	
	\par Now as there are only two distinct zero sets, there cannot exist any cycle in the graph $\Gamma(C(X)_F)$. Hence the girth $gr(\Gamma(C(X)_F))=\infty$.
	\par Now let $diam(\Gamma(C(X)_F))=2$ or girth $gr(\Gamma(C(X)_F))=\infty$. By the above theorem, we see that if $X$ contains more than two points then diameter and girth both are $3$. Hence we have $|X|=2$.
\end{proof}
\begin{defn}
	For a vertex $f$ in a graph $\Gamma(C(X)_F)$ the associated number $e(f)$ is defined by $e(f)=\max\{d(f,g):g\neq f\}$. The vertex $g$ with smallest associated number is called a centre of the graph. The  associated number of the centre vertex is called the radius of the graph and it is denoted by $\rho(\Gamma(C(X)_F))$. 
\end{defn}
\par The next result is about the associated number of any vertex in the graph $\Gamma(C(X)_F)$.
\begin{lem}
	For any vertex $f$ in the graph $\Gamma(C(X)_F)$, we have 
	\[e(f)=\left\{\begin{array}{ll}
		2$\text{ if } $X\setminus \mathcal{Z}(f)\text{ is singleton}\\
		3 \text{ otherwise.}
	\end{array}\right.\]
\end{lem}
\begin{proof}
	Let $X\setminus\mathcal{Z}(f)=\{x_0\}$. Now consider $g\in \Gamma(C(X)_F)$ such that $g\neq f$. Then there are only two possibilities, namely $\mathcal{Z}(g)$ contains $x_0$ or does not contain $x_0$. If $\mathcal{Z}(g)$ contains $x_0$ then $fg=0$. In this case $f$ and $g$ are adjacent to each other. Hence $d(f,g)=1$. If $\mathcal{Z}(g)$ does not contain $x_0$ then $\mathcal{Z}(g)\subseteq \mathcal{Z}(f)$. Then $\mathcal{Z}(f)\cap \mathcal{Z}(g) = \mathcal{Z}(g)\neq \emptyset$ and $\mathcal{Z}(f)\cup \mathcal{Z}(g)=\mathcal{Z}(f)\neq \emptyset$. Thus by Theorem \ref{deg-thm}, $d(f,g)=2.$ So we conclude that $e(f)=2.$
	\par On the other hand, let $X\setminus\mathcal{Z}(f)$ contains at least two points, namely $x_0$ and $y_0$. By Theorem \ref{deg-thm}, we see that $e(f)\leq 3.$ Now choose $g=1-\chi_{\{x_0\}}$. Then $\mathcal{Z}(g)=\{x_0\}$. Clearly, $\mathcal{Z}(f)\cap \mathcal{Z}(g)=\emptyset$ and $\mathcal{Z}(f)\cup \mathcal{Z}(g)\neq \emptyset$ because $y_0$ does not belong to the union. So by Theorem \ref{deg-thm}, for this particular $g$, we get $d(f,g)=3$. Hence we conclude that $e(f)=3$.
\end{proof}
\begin{cor}\label{radius}
	The radius of the graph $\Gamma(C(X)_F)=2$, that is, $\rho(\Gamma(C(X)_F))=2.$
\end{cor}
\begin{proof}
	We can always consider a function $f=\chi_{\{x_0\}}$ then $X\setminus \mathcal{Z}(f)=\{x_0\},$ that is, singleton. Hence $e(f)=2$. So we have radius of $\Gamma(C(X)_F)=2$.
\end{proof}
%\subsection{Cycle structure in $C(X)_F$}
\begin{defn}
	A graph $G$ is said to be \\
	$(i)$ triangulated if every vertex of the graph $G$ is a vertex of a triangle.\\	
$(ii)$ hyper-triangulated if every  edge of the graph $G$ is an edge of a triangle.
\end{defn}
\begin{thm}\label{try}
	The graph $\Gamma(C(X)_F)$ is neither triangulated nor hyper-triangulated.
\end{thm}
\begin{proof}
	Let us consider a point $x_0\in X$. Now consider two functions $f=\chi_{x_0}$ and $g=1-\chi_{x_0}$. Then $\mathcal{Z}(f)\cup\mathcal{Z}(g)=X$ and $\mathcal{Z}(f)\cap \mathcal{Z}(g)=\emptyset$. Then by Lemma \ref{deg 2} it is not possible to get a cycle of length 3 that contains $f$ and the edge connecting $f$ and $g$. So the graph $\Gamma(C(X)_F)$ is neither triangulated nor hyper-triangulated.
\end{proof}
The above result is significantly different in the case of $C(X)$. In fact, we have 
\begin{proposition}[\cite{zero}]
	\begin{enumerate}[$(i)$]
		\item $\Gamma(C(X))$ is triangulated if and only if $X$ does not contain any non-isolated points. 
		\item $\Gamma(C(X))$ is hyper-triangulated if and only if $X$ is a connected middle $P$-space.
	\end{enumerate}
\end{proposition}
 For definition of middle $P$-space see \cite{zero}.
\begin{defn}
	For two vertices $f$ and $g$ in the graph $\Gamma(C(X)_F)$, we denote by $c(f,g)$ the length of the smallest cycle containing $f$ and $g$. If there is no cycle containing $f$ and $g$, we declare $c(f,g)=\infty$.
\end{defn}
\par In the following theorem, we shall discuss all possible values of  $c(f,g)$ in the graph $\Gamma(C(X)_F)$.
\begin{thm}
	Let $f$ and $g$ be two vertices in the graph $\Gamma(C(X)_F)$. Then \begin{enumerate}[$(i)$]
		\item $c(f,g)=3$ if and only if $\mathcal{Z}(f)\cup\mathcal{Z}(g)=X$ and $\mathcal{Z}(f) \cap \mathcal{Z}(g)\neq \emptyset$.
		\item $c(f,g)=4$ if and only if $\mathcal{Z}(f)\cup\mathcal{Z}(g)=X$ and $\mathcal{Z}(f) \cap \mathcal{Z}(g)= \emptyset$ or $\mathcal{Z}(f)\cup\mathcal{Z}(g)\neq X$ and $\mathcal{Z}(f) \cap \mathcal{Z}(g) \neq \emptyset$.
		\item $c(f,g)=6$ if and only if $\mathcal{Z}(f)\cup\mathcal{Z}(g)\neq X$ and $\mathcal{Z}(f) \cap \mathcal{Z}(g) = \emptyset$.
	\end{enumerate}
\end{thm} 
\begin{proof}
		(i) Let $\mathcal{Z}(f)\cup\mathcal{Z}(g)=X$ and $\mathcal{Z}(f) \cap \mathcal{Z}(g)\neq \emptyset$. Then by Lemma \ref{deg1} and Lemma \ref{deg 2} we can say that $f$ and $g$ are adjacent to each other and there is another function $h$ adjacent to both $f$ and $g$. So we get a triangle with vertices $f,g$ and $h$. This implies that $c(f,g)=3$. Conversely, if $c(f,g)=3$ then there is a triangle having $f,g$ and $h$ as its vertices for some other vertex $h$. Now appealing to Lemma \ref{deg1} and Lemma \ref{deg 2}, we can say that $\mathcal{Z}(f)\cup\mathcal{Z}(g)=X$ and $\mathcal{Z}(f) \cap \mathcal{Z}(g)\neq \emptyset$.\\
		
		(ii) Let $\mathcal{Z}(f)\cup\mathcal{Z}(g)=X$ and $\mathcal{Z}(f) \cap \mathcal{Z}(g)= \emptyset$. Then by Lemma \ref{deg1}, $f$ and $g$ are adjacent to each other. Now by Lemma \ref{deg 3}, there are vertices $h_1$ and $h_2$ such that $f$ is adjacent to $h_1$, $h_1$ is adjacent to  $h_2$ and $h_2$ is adjacent to $g$. So we get a quadrilateral with vertices in order, $f,h_1,h_2$ and $ g$. As $\mathcal{Z}(f) \cap \mathcal{Z}(g)= \emptyset$, by Lemma \ref{deg 2}, there is no triangle containing $f$ and $g$ as its vertices. So $c(f,g)=4$.
		
		\par Now let $\mathcal{Z}(f)\cup\mathcal{Z}(g)\neq X$ and $\mathcal{Z}(f) \cap \mathcal{Z}(g) \neq \emptyset$. Then by Lemma \ref{deg1}, $f$ and $g$ are not adjacent to each other. By Lemma \ref{deg 2}, there is a vertex $h$ such that $h$ is adjacent to both $f$ and $g$. Clearly then $2h$ is also adjacent to both $f$ and $g$. So we get a quadrilateral containing vertices in order $f,h,g$ and $2h$. Again condition $\mathcal{Z}(f)\cup\mathcal{Z}(g)\neq X$ implies that it is not possible to have a triangle containing $f$ and $g$ as its vertices. So $c(f,g)=4.$
		\par To prove the converse, let $c(f,g)=4.$ Now $\mathcal{Z}(f)\cup \mathcal{Z}(g)=X$, then we must have  $\mathcal{Z}(f)\cap \mathcal{Z}(g)=\emptyset$, otherwise we have a triangle having vertices $f$ and $g$. If we have $\mathcal{Z}(f)\cup \mathcal{Z}(g)\neq X$, then $f$ and $g$ are not adjacent to each other. But there is a quadrilateral containing $f$ and $g$. So there must exist two functions $h_1$ and $h_2$ such that both $h_1$ and $h_2$  are adjacent to both $f$ and $g$. So by Lemma \ref{deg 2}, we have $\mathcal{Z}(f) \cap \mathcal{Z}(g) \neq \emptyset$.\\
		
		(iii) Let $\mathcal{Z}(f)\cup\mathcal{Z}(g)\neq X$ and $\mathcal{Z}(f) \cap \mathcal{Z}(g) = \emptyset$. Then $f$ and $g$ are not adjacent to each other. As $
		\mathcal{Z}(f) \cap \mathcal{Z}(g) = \emptyset$, by Lemma \ref{deg 3}, there are two vertices $h_1$ and $h_2$ in $\Gamma (C(X)_F)$ such that there is a path connecting $f,h_1,h_2 $ and $g$ in order. So immediately there is another path connecting $g,2h_2,2h_1$ and $f$. So we get a cycle of length 6, namely $f,h_1,h_2,g,2h_2,2h_1$ and $f$. Let us make it clear that with the given condition it is not possible to get a cycle of length 5. As $f$ and $g$ are not adjacent to each other, to have a cycle of length 5, we must have a path of length 2 joining $f$ and  $g$ which is not possible as $\mathcal{Z}(f)\cap\mathcal{Z}(g)=\emptyset.$ This implies that $c(f,g)=6.$
		\par Conversely, let $c(f,g)=6.$ Then by proof of (i) and (ii), we have $\mathcal{Z}(f)\cup\mathcal{Z}(g)\neq X$ and $\mathcal{Z}(f) \cap \mathcal{Z}(g) = \emptyset$. 
\end{proof}
\section{Conclusion} In conclusion, we point out the resemblance and difference between the behaviour of the rings $C(X)_F$ and $C(X)$, based on the results obtained in the article. Though $B_1(X)$ is always an overring of $C(X)$, we have observed in this paper that $B_1(X)$ is not always an overring of $C(X)_F$. We have proved that in case $X$ is perfectly normal, $B_1(X)$ is indeed an overring of $C(X)_F$ [Theorem \ref{perfect}]. As opposed to the fact that $C(X)$ and $B_1(X)$ both are closed under uniform limits, we have observed [Theorem \ref{uniform}] that $C(X)_F$ is closed under uniform limits whenever $X$ has only finitely many non-isolated points. There are quite a lot of analogies between the zero divisor graphs of $C(X)$ and $C(X)_F$. However, there are a few deviations in results obtained for the zero divisor graph of $C(X)_F$ from those of the zero divisor graph of $C(X)$. For example, it is known that $\Gamma(C(X))$ is triangulated if $X$ has no isolated point and it is hyper-triangulated if $X$ is connected and a middle $P$-space \cite{zero}. But in case of $\Gamma(C(X)_F)$, we  have established that $\Gamma(C(X)_F)$  is never triangulated and never hyper-triangulated. Also, the radius of $\Gamma(C(X)_F)$ is always $2$ irrespective of the nature of $X$, whereas the radius  of $\Gamma(C(X))$ takes the value  $3$ when $X$ has no isolated point and 2 otherwise.
% [Corollary \ref{radius}].


\begin{thebibliography}{6}
\bibitem{blum} M. R. Ahmadi Zand, \textit{An algebraic characterization of Blumberg spaces}, Quaest. Math. 33 (2) (2010), 1-8.
\bibitem{appl} M. R. Ahmadi Zand and Z. Khosravi,  \textit{Remarks on the rings of functions which have a finite number of discontinuities}. Appl. Gen. Topol. 22(1) (2021), 139-147.
 \bibitem{fdapsl} D. F. Anderson and P. S. Livingston, \textit{The zero-divisor graph of commutative ring}, J. Algebra, 217 (1999), 434–447.
\bibitem{zero} F. Azarpanah, and M. Motamedi, \textit{Zero-divisor graph of C(X)}, Acta Math. Hungar., 108 (1-2) (2005), 25-36. 
\bibitem{bare1} A. Deb Ray and A. Mondal, \textit{On rings of Baire one functions}, Appl. Gen. Topol. 20 (1) (2019), 237–249.
\bibitem{bare2} A. Deb Ray and A. Mandal, \textit{Ideals in $B_1(X)$ and residue class rings of $B_1(X)$ modulo an ideal}, Appl. Gen. Topol. 20 (2) (2019), 379-393.
\bibitem{gill} L. Gillman and M. Jerison, \textit{Rings of Continuous Functions}, Springer, London. 1976.
\bibitem{main} Z. Gharebaghi, M. Ghirati and A. Taherifar,  \textit{On the ring of functions which are discontinuous on a finite set}, Houston J. Math. 44 (2) (2018), 721–739.
\bibitem{vesel} L. Vesely, \textit{Characterization of Baire-one functions between topological spaces}, Acta Univ. Carolin. Math. Phys., 33 (2) (1992), 143–156.

\end{thebibliography}
\end{document}